\documentclass[a4paper]{amsart}

\usepackage{amsmath,amsthm,amsfonts,amssymb,amscd}
\input{xy}
\xyoption{all}

\newtheorem{theorem}{Theorem}[section]

\newtheorem{defi}[theorem]{Definition}
\newtheorem{rem}[theorem]{Remark}
\newtheorem{lem}[theorem]{Lemma}

\newcommand{\dem}{{\textit{Proof. }}}

\newenvironment{sis}{\left\{\begin{aligned}}{\end{aligned}\right.}

\numberwithin{equation}{section}

\newcommand{\Z}{\mathbb{Z}}
\newcommand{\N}{\mathbb{N}}

\newcommand{\0}{\underline{0}}

\newcommand{\g}{\mathfrak{g}}

\newcommand{\ep}{\epsilon}

\renewcommand{\d}{{\rm d}}

\newcommand{\Sq}{{\rm Sq}}

\newcommand{\ti}{\widetilde}

\begin{document}

\title[Restricted infinit. deform. of restricted simple Lie algebras]{Restricted Infinitesimal Deformations of Restricted Simple Lie Algebras}

\author{Filippo Viviani}


\address{Dipartimento di Matematica,
Universit\`a Roma Tre,
Largo S. Leonardo Murialdo 1,
00146 Roma (Italy)}
\email{viviani@mat.uniroma3.it}



\begin{abstract}
We compute the restricted infinitesimal deformations of the restricted simple Lie algebras
over an algebraically closed field of characteristic $p\geq 5$.
\end{abstract}

\maketitle

\section{Introduction}

\emph{Simple Lie algebras} over an algebraically closed field
$F$ of characteristic $p\neq 2,3$ have recently been classified by
Block-Wilson-Premet-Strade (see  \cite{BW1}, \cite{SW}, \cite{STR1},\cite{STR2},
\cite{STR3},\cite{STR4},\cite{STR5},\cite{STR6},\cite{PS1}, \cite{PS2}, \cite{PS3},
\cite{PS4}, \cite{PS5}, \cite{PS6}, \cite{PS7}, \cite{STR}).
The classification says that for $p\geq 7$ the simple Lie algebras
can be of two types: of classical type and of generalized Cartan type.

\hspace{0,4cm}
The algebras of \emph{classical type} are obtained by considering the simple Lie
algebras in characteristic zero (classified via Dynkin diagrams), by taking a model
over the integers via the Chevalley bases and then reducing modulo the prime $p$
(see \cite{SEL}).

\hspace{0,4cm}
The algebras of \emph{generalized Cartan type} were constructed by
Kostrikin-Shafare\-vich, Wilson and Kac (\cite{KS}, \cite{KS2}, \cite{WIL1}, 
\cite{KAC1}, \cite{KAC2}, \cite{KAC4}, \cite{KAC3}, \cite{WIL2})
and are divided into four families, called Witt-Jacobson, special,
Hamiltonian and contact algebras. These four families are the finite-dimen\-sional
analogue of the four classes of infinite-dimensional complex simple Lie algebras,
which occurred in Cartan's classification (\cite{CAR})
of Lie pseudogroups.

\hspace{0,4cm}
In characteristic $p=5$, apart from the above two types of algebras, there is one
more family of simple Lie algebras called \emph{Melikian algebras} (introduced in \cite{MEL}).
In characteristic $p=2, 3$, there are many exceptional simple Lie algebras
(see \cite[page 209]{STR}) and the classification seems still far away.

\hspace{0,4cm}
We are interested in a particular class of modular Lie algebras
called \emph{restricted}. These can be characterized as those
modular Lie algebras such that the $p$-power of an inner derivation
(which in characteristic $p$ is a derivation) is still inner.
Important examples of restricted Lie algebras are the ones coming
from groups schemes. Indeed there is a one-to-one correspondence
between restricted Lie algebras and finite group schemes whose
Frobenius vanishes (see \cite[Chap. 2]{DG}).


\hspace{0,4cm}
The aim of this paper is to compute the \emph{restricted infinitesimal deformations}
of the restricted simple Lie algebras in characteristic $p\geq 5$.
By standard facts of deformation theory (see for example \cite{GER1}),
restricted infinitesimal deformations of a
restricted Lie algebra $\g$ are parametrized by
the second restricted cohomology group $H_*^2(\g,\g)$ of $\g$ with values in the adjoint
representation (see \cite{HOC}).



\hspace{0,4cm}
The restricted simple Lie algebras of classical type are known to be rigid
as Lie algebras, under the assumption $p\geq 5$ (see \cite{RUD}). This is equivalent
to the vanishing of the second ordinary cohomology group
$H^2(\g, \g)$ for $\g$ restricted simple of classical type.
Therefore, using the so-called Hochschild $6$-term exact sequence
(see (\ref{6-term}) below), one can easily deduce the vanishing of $H_*^2(\g,\g)$, which
implies that these algebras are rigid also as restricted Lie algebras.
Note that in characteristic zero, the vanishing of $H^2(\g, \g)$
for $\g$ simple follows from a  classical result, known as second
Whitehead's Lemma (see for example \cite{HiSt}). Interestingly, some of the
classical simple Lie algebras admit non-trivial deformations in characteristic $p=2$ or $3$
(see \cite{DZU1}, \cite{Che2}, \cite{Che3}, \cite{Che1}).


\hspace{0,4cm}
Under the assumption $p\geq 5$, we compute the restricted infinitesimal
deformations of the restricted simple Lie algebras not of classical type: the four infinite
families $W(n):=W(n,\underline{1})$, $S(n):=S(n,\underline{1})$, $K(n):=K(n,\underline{1})$,
$H(n):=H(n,\underline{1})$ and the exceptional restricted Melikian algebra
$M:=M(1,1)$
in characteristic $p=5$.
Using the notations about those algebras which we are going to recall in what follows
and the squaring operation ${\rm Sq}$ (see section 2.2),
we can state our results as follows.

\begin{theorem}\label{W-finaltheorem}
The infinitesimal restricted deformations of the restricted Ja\-cobson-Witt algebra
$W(n)$ are given by
$$H_{*}^2(W(n),W(n))=H^2(W(n),W(n))=\bigoplus_{i=1}^n \langle {\rm Sq}(D_i)\rangle_F.$$
\end{theorem}

\begin{theorem}\label{S-finaltheorem}
The infinitesimal restricted deformations of the restricted special algebra $S(n)$ are given by
$$H_*^2(S(n),S(n))=\bigoplus_{i=1}^n \langle {\rm Sq}(D_i)\rangle_F.$$
\end{theorem}

\begin{theorem}\label{K-finaltheorem}
The infinitesimal restricted deformations of the restricted contact algebra $K(n)$ are given by
$$H_*^2(K(n),K(n))=H^2(K(n),K(n))=\bigoplus_{i=1}^{2m} \langle {\rm Sq}(x_i) \rangle_F \oplus
\langle {\rm Sq}(1)\rangle_F.$$
\end{theorem}

\begin{theorem}\label{H-finaltheorem}
The infinitesimal restricted deformations of the restricted Ha\-miltonian algebra
$H(n)$ are given by
$$H_*^2(H(n),H(n))=\bigoplus_{i=1}^n \langle {\rm Sq}(x_i)\rangle_F\oplus
\langle \Phi \rangle_F,  $$
where the cocycle $\Phi$ is defined as
$$\Phi(x^a,x^b)=\sum_{\stackrel{\0<\delta\leq a, \widehat{b}}{|\delta|=3}}
\binom{a}{\delta}\binom{b}{\widehat{\delta}}\sigma(\delta)\:\delta!\:
x^{a+b-\delta -\widehat{\delta}}.
$$
\end{theorem}

\begin{theorem}\label{M-finaltheorem}
The restricted infinitesimal deformations of the restricted Melikian algebra $M$ are given by
$$H_{*}^2(M,M)=H^2(M,M)=\langle {\rm Sq}(1)\rangle_F \bigoplus_{i=1}^2 \langle {\rm Sq}(D_i)
\rangle_F \bigoplus_{i=1}^2 \langle {\rm Sq}(\ti{D_i})\rangle_F.$$
\end{theorem}

We prove these theorems using
the computation of the second ordinary cohomology group $H^2(\g,\g)$ that we
performed in \cite{VIV1}, \cite{VIV2} and \cite{VIV3}, together with
the $6$-term exact sequence of Hochschild (see section 2.1) that relates the
ordinary and restricted cohomology.

\hspace{0,4cm}
Note that the second restricted cohomology group for the above algebras
is freely generated over $F$ by the squaring operators of the elements of negative
degree, with one remarkable exception: for the Hamiltonian algebra $H(n)$, there is
an exceptional extra-cocycle $\Phi$.
Observe moreover that, quite interestingly, for the special algebras $S(n)$ and the
Hamiltonian algebras $H(n)$, the restricted cohomology group $H_*^2(\g,\g)$ is a
proper subgroup of the ordinary cohomology group $H^2(\g,\g)$. In other words,
there are infinitesimal deformations of the algebra that do not admit a restricted
structure.

\hspace{0,4cm}
Since the restricted infinitesimal deformations of
a restricted Lie algebra $\g$ correspond to the infinitesimal deformations
of the associated finite group scheme $G$ of height one, the above results give
the infinitesimal deformations of some simple finite group schemes
(see \cite{VIV4} and \cite{VIV5} for more details).
In order to complete the picture, it would be
very interesting to extend the above computations to the minimal
$p$-envelope of all the simple Lie algebras.

\hspace{0,4cm}
The paper is organized as follows.
The section $2$ contains preliminary results. We quickly review the ordinary
and restricted cohomology of Lie algebras and the $6$-term
Hochschild exact sequence relating the first two ordinary and restricted cohomology
groups. Moreover we recall the definition of the squaring operation.
In each of the remaining five sections, we recall the basic definitions of the
five classes of restricted simple Lie algebras of non-classical type and we compute the
corresponding infinitesimal deformations.

\section{Preliminaries}

\subsection{Ordinary and restricted cohomology}

In this section we review, in order to fix notations, the ordinary and restricted
cohomology of Lie algebras, following \cite{HS} and \cite{HOC}.

\hspace{0,4cm}
Let $\g$ be a Lie algebra over a field $F$. We denote by $U_{\g}$ the universal
enveloping algebra
of $\g$ and by $I_{\g}$ its augmentation ideal, that is the kernel of the augmentation map
$\ep:U_{\g}\to F$. For a $\g$-module $M$ (or, equivalently, an unital $U_{\g}$-module),
the (ordinary) cohomology groups $H^n(\g,M)$ are the right derived functor of the fixed point
functor $M\mapsto M^{\g}$, considered as a functor from the category of $\g$-modules
to the category of abelian groups.
They can be computed into two different ways, using a Lie complex or an associative complex.

\hspace{0,4cm}
The Lie complex has $n$-dimensional cochains $C^n(\g, M)=\{f: \Lambda^{n}(\g)\to M\}$
and differential $\d:C^n(\g,M)\to C^{n+1}(\g,M)$ defined by
\begin{equation*} \begin{split}
\d f(x_0,\dots,x_n)=&\sum_{i=0}^n(-1)^i x_i\cdot f(x_0,\dots, \hat{x_i},
\dots, x_n)+\\
&\sum_{p<q}(-1)^{p+q}f([x_p,x_q],x_0,\dots,\hat{x_p},\dots, \hat{x_q},
\dots x_n)
\end{split} \end{equation*}
where the sign $\hat{}$ means that the argument below must be omitted.

\hspace{0,4cm}
The associative complex has $n$-dimensional cochains $C^n(I_{\g},M)=\{g:I_{\g}^{\otimes n}
\to M\}$, and differential $\d:C^n(I_{\g},M)\to C^{n+1}(I_{\g},M)$ defined by
\begin{equation*}
\d g(s_0,\cdots s_n)=s_0\cdot g(s_1,\cdots, s_n)+\sum_{i=1}^n (-1)^i g(s_0,\cdots,s_{i-1}s_i,
\cdots, s_n).
\end{equation*}

\hspace{0,4cm}
Now let $(\g,[p])$ be a restricted Lie algebra over $F$.
Denote by $U_{\g}^{[p]}:=U_{\g}/(x^p-x^{[p]})$ the restricted enveloping algebra of
$(\g, [p])$ and with $I_{\g}^{[p]}$ its augmented ideal. For a restricted $\g$-module $M$
(or, equivalently, an unital $U_{\g}^{[p]}$-module),
the restricted cohomology groups $H_{*}^n(\g,M)$ are the right derived functor
of the fixed point
functor $M\mapsto M^{\g}$, considered as a functor from the category of restricted $\g$-modules
to the category of abelian groups. Explicitly, these can be calculated via an associative
complex which is obtained from the one described above for ordinary cohomology groups simply
by replacing $I_{g}$ with $I_{\g}^{[p]}$. Observe, on the other hand, that the ordinary Lie
complex does not generalize to restricted cohomology, a fact which makes the computation of
the restricted cohomology harder than the ordinary one.

\hspace{0,4cm}
There is a $6$-term exact sequence relating the first two ordinary and restricted
cohomology groups (see \cite{HOC}):
\begin{equation}\label{6-term}
\begin{split}
&0\to H_{*}^1(\g,M)\to H^1(\g,M)\stackrel{D}{\longrightarrow}
{\rm Hom}_{Fr}(\g, M^{L})\to\\
& \to H_{*}^2(\g,M) \to H^2(\g,M)\stackrel{H}{\longrightarrow}
{\rm Hom}_{Fr}(\g,H^1(\g,M)),
\end{split}
\end{equation}
where ${\rm Hom}_{Fr}(V,W)$ denotes the Frobenius-semilinear morphisms between
the two $F$-vector spaces $V$ and $W$, that is
$${\rm Hom}_{Fr}(V,W)=\{f:V\to W\: | \: f(\alpha x+\beta y)= \alpha^p
f(x)+\beta^p f(y) \},$$
for any $ \alpha, \beta \in F, x,y\in V$
and $D$ and $H$ are defined on the Lie cochains $\phi\in H^1(\g,M)$ and
$\psi\in H^2(\g,M)$ as, respectively (for any $x, y\in \g$):
\begin{equation*}\label{DH-morphisms}
\begin{sis}
& D_{\phi}(x)=x^{p-1}\circ \phi(x)-\phi(x^{[p]}), \\
& H_{\psi}(x)\cdot y=\sum_{j=0}^{p-1} x^j \circ \psi(x, ({\rm ad}\,x)^{p-1-j}(y))-
\psi(x^{[p]},y).
\end{sis}
\end{equation*}
In the particular case in which $M=\g$ is the adjoint representation and the
algebra $\g$ has no center, the above $6$-term exact sequence (\ref{6-term}) becomes
\begin{equation}\label{inclusion-sequence}
\begin{sis}
& H_{*}^1(\g,\g)=H^1(\g,\g),\\
& 0\to H_{*}^2(\g,\g)\to H^2(\g,\g)\stackrel{H}{\longrightarrow}
{\rm Hom}_{Fr}(\g,H^1(\g,\g)),
\end{sis}
\end{equation}
and the operator $H$ becomes (for any $\psi\in H^2(\g,\g)$ and $x, y\in \g$)
\begin{equation}\label{H-morphism}
 H_{\psi}(x)\cdot y=\sum_{j=0}^{p-1} ({\rm ad}\,x)^j \circ \psi(x, ({\rm ad}\,x)^{p-1-j}(y))-
\psi(x^{[p]},y).
\end{equation}

\begin{rem}
The Hochschild $6$-term exact sequence (\ref{6-term})
has been interpreted as the initial sequence
of two different spectral sequences relating the ordinary and restricted cohomology:
$$\begin{sis}
&E_1^{p,q}={\rm Hom}_{Fr}(S^p \g, H^{q-p}(\g,M))\Rightarrow H_{*}^{p+q}(\g,M) \:
\text{ if } p\neq 2 \hspace{0,5cm} & \text{\cite{JAN}},\\
&E_2^{p,q}={\rm Hom}_{Fr}(\Lambda^q \g, H_{*}^{p}(\g,M))\Rightarrow H^{p+q}(\g,M)
\hspace{1cm} & \text{\cite{FAR}},
\end{sis}$$
where $S^p\g$ and $\Lambda^q\g$ denote, respectively, the $p$-th symmetric power and
the $q$-th alternating power.
\end{rem}

\subsection{Squaring operation}

There is a canonical way to produce $2$-cocycles in $H^2(\g,\g)$ over a field of characteristic
$p>0$, namely the squaring operation (see \cite{GER1}).
Given a derivation $\gamma$ (inner or not), one defines the squaring of
$\gamma$ to be
\begin{equation}\label{Square}
\Sq(\gamma)(x,y)=\sum_{i=1}^{p-1}\frac{[\gamma^i(x),\gamma^{p-i}(y)]}{i!(p-i)!}
\end{equation}
where $\gamma^i$ is the $i$-iteration of $\gamma$.
In \cite{GER1} it is shown that $[\Sq(\gamma)]\in H^2(\g,\g)$ is an obstruction to integrability
of the derivation $\gamma$, that is to the possibility of finding an automorphism
of $\g$ extending the infinitesimal automorphism given by $\gamma$.

\section{The Witt-Jacobson algebra}

Let us recall the definition of the restricted Witt-Jacobson algebra, following
\cite[Section 4.2]{FS}.
Let $A(n)=A(n;\underline{1}):=F[x_1,\dots , x_n]/(x_1^p, \dots, x_n^p)$ be the ring of
$p$-truncated polynomials
in $n$ variables over a field $F$ of positive characteristic $p\geq 5$.

\begin{defi}
The restricted Witt-Jacobson algebra $W(n)=W(n;\underline{1})$ is the restricted Lie algebra
${\rm Der}_FA(n)$ of derivations of $A(n)=F[x_1,\ldots,x_n]/$ $(x_1^p,$ $\ldots,$ $x_n^p)$.
\end{defi}

The Witt-Jacobson algebra $W(n)$ is
a free $A(n)$-module with basis $\{D_1,$ $ \dots, $ $D_n\}$, where we put
$D_j:=\frac{\partial}{\partial x_j}$. Therefore ${\rm dim}_F(W(n))=np^n$ with a
basis over $F$ given by $\{x^aD_j\: | \: 1\leq j\leq n, \: 0\leq a_i \leq p-1\}$.
It has a natural grading obtained by assigning to the element $x^aD_j$ the degree
$|a|-1:=\sum_{i=1}^n a_i-1$. In particular the elements of negative degree are
$W(n)_{-1}=\langle D_1,\cdots, D_n\rangle_F$.
The $[p]$-map is defined on the elements of the base by
$$(x^aD_j)^{[p]}=\begin{sis}
&x^aD_j & \text{ if } \: x^aD_j=x_jD_j, \\
&0 & \text{ otherwise. } \\
\end{sis}$$
It is a classical result of Celousov (see \cite{CEL} or \cite[Section 4.8]{FS}) that
every derivation of $W(n)$ is inner or in other words that
\begin{equation}
H_{*}^1(W(n),W(n))=H^1(W(n),W(n))=0.
\end{equation}
Therefore, from the Hochschild exact sequence (\ref{inclusion-sequence}), we deduce
that $H_*^2(W(n),$ $W(n))=H^2(W(n),W(n))$. The Theorem \ref{W-finaltheorem}
follows from \cite[Theorem 1.1]{VIV1}:
$$H^2(W(n),W(n))=\bigoplus_{i=1}^n \langle {\rm Sq}(D_i)\rangle_F.$$

\section{The special algebra}

Let us recall the definition of the restricted special algebra, following
\cite[Section 4.3]{FS}.
Fix an integer $n\geq 3$ and a field $F$ of characteristic $p\geq 5$.
Consider the following map, called divergence:

$${\rm div}:\begin{sis}
 W(n) & \to A(n)\\
 \sum_{i=1}^n f_iD_i & \mapsto \sum_{i=1}^n D_i(f_i).
\end{sis}$$
The kernel of the divergence map
$S'(n)=S'(n;\underline{1})=\{E\in W(n)\: | \: {\rm div}(E)=0\}$ is a graded subalgebra of $W(n)$
of dimension $(n-1)p^n+1$.

\begin{defi}
The restricted special algebra $S(n)=S(n;\underline{1})$ is the derived
subalgebra of $S'(n)$:
$$S(n):=S'(n)^{(1)}=[S'(n),S'(n)].$$
\end{defi}

It turns out that there is an exact sequence
\begin{equation*}\label{S-S'-sequence}
0\to S(n)\to S'(n)\to \oplus_{i=1}^n \langle x^{\tau-(p-1)\epsilon_i}D_i
\rangle_F \to 0,
\end{equation*}
where $\tau:=(p-1,\cdots,p-1)$ and $\ep_i$ is the $n$-tuple having $1$ at the
$i$-th place and $0$ outside. Therefore $S(n)$ has $F$-dimension
$(n-1)(p^n-1)$. A set of generators (but not linearly independent!) of $S(n)$ is given
by the elements $\{D_{ij}(f)\: | \: f\in A(n), \: 1\leq i <j\leq n\}$,
where the maps $D_{ij}$ are defined by:
$$D_{ij}:\begin{sis}
A(n)& \longrightarrow W(n)\\
f & \mapsto D_j(f)D_i-D_i(f)D_j.\\
\end{sis}$$
In particular, the elements of negative degree are $S(n)_{-1}=
\langle D_1,\cdots, D_n\rangle_F$. The $[p]$-map on the above generators is given by
$$D_{ij}(x^a)^{[p]}=\begin{sis}
& D_{ij}(x^a) & \text{ if } x^a=x_ix_j, \\
& 0 & \text{ otherwise. }\\
\end{sis}$$
The first cohomology group of the adjoint representation is equal to (see \cite{CEL}
or \cite[Section 4.8]{FS}):
\begin{equation*}
H_{*}^1(S(n),S(n))=H^1(S(n),S(n))=\bigoplus_{i=1}^n
\langle {\rm ad}(x^{\tau-(p-1)\epsilon_i}D_i)\rangle_F \bigoplus
\langle {\rm ad}(x_1D_1)\rangle_F.
\end{equation*}
From this result, we can deduce a criterion saying when a derivation  of $S(n)$
is inner. First we introduce the following notation. Observe that,
expressing any element $E\in S(n)$ as $F$-linear combination of the generators
$D_{ij}(x^{\tau})$,
the coefficients of the terms of minimal degree $D_{ij}(x_j)=D_i$ and of maximal degree
$D_{ij}(x^{\tau})=x^{\tau-\ep_j}D_i-x^{\tau-\ep_i}D_j$ are well-defined,
that is they are the same for any such expression of $E$.
We call the above coefficients  $E_{D_i}$ and $E_{D_{ij}(x^{\tau})}$, respectively.

\begin{lem}\label{S-inner}
A derivation $\gamma:S(n)\to S(n)$ is inner if and only it satisfies the
following conditions:
\begin{itemize}
\item[(i)] For every $1\leq i\leq n$, there exists a $j\neq i$ such that
$\gamma(x_i^{p-1}D_j)_{D_{ij}(x^{\tau})}=0.$
\item[(ii)] $\sum_{k=1}^n \gamma(D_k)_{D_k}=0$.
\end{itemize}
\end{lem}
\dem
We first prove that the two conditions are necessary.
Consider an inner derivation ${\rm ad}(D)$, with $D\in S(n)\subset W(n)$.
To prove condition (i), write $D$ as linear combination of the base elements $x^aD_h\in W(n)$.
Consider the element
of $W(n)$ given by (for $i\neq j$)
$${\rm ad}(x^aD_h)(x_i^{p-1}D_j)=[x^aD_h,x_i^{p-1}D_j]=x^aD_h(x_i^{p-1})D_j-
x_i^{p-1}D_j(x^a)D_h.$$
Clearly, the two elements at the end cannot be equal to $x^{\tau-\ep_j}D_i$
and therefore ${\rm ad}(D)(x_i^{p-1}D_j)_{D_{ij}(x^{\tau})}=0$.
To prove the condition (ii), write $D=\sum_{i=1}^n a_i x_iD_i+E$ with
$E_{x_iD_i}=0$ for every $i$.
Clearly $D\in S(n)$ if and only if $E\in S(n)$ and
$\sum_{i=1}^n a_i=0$. We compute
$$\sum_{k=1}^n {\rm ad}(D)(D_k)_{D_k}=-\sum_{k=1}^na_k=0.$$
However, the two conditions are also sufficient since
$$\begin{sis}
&\sum_{k=1}^n{\rm ad}(x_1D_1)(D_k)_{D_k}={\rm ad}(x_1D_1)(D_1)_{D_1}=-1,\\
&{\rm ad}(x^{\tau-(p-1)\epsilon_i}D_i)(x_i^{p-1}D_j)=-x^{\tau-\ep_i}D_j+x^{\tau-\ep_j}D_i=
-D_{ij}(x^{\tau}).\qed \\
\end{sis} $$

In \cite[Theorem 1.2]{VIV1}, we prove that the second ordinary cohomology
group of the adjoint representation of $S(n)$ is
\begin{equation}\label{S-ord-def}
H^2(S(n),S(n))=\bigoplus_{i=1}^n \langle {\rm Sq}(D_i)\rangle_F   \bigoplus
\langle\Theta \rangle_F,
\end{equation}
where $\Theta$ is defined by $\Theta(D_i,D_j)=D_{ij}(x^{\tau})$ and extended by $0$
elsewhere.
Using this result and the Hochschild exact sequence (\ref{inclusion-sequence}), we can compute
the second restricted cohomology group.

\dem[of Theorem \ref{S-finaltheorem}]
The cocycle $\Theta$ does not belong to $H_{*}^2(S(n),S(n))$. Indeed, using that
$D_i^{[p]}=0$, we compute (for $i\neq j$)
$$H_{\Theta}(D_i)\cdot x_i^{p-1}D_j=\sum_{k=0}^{p-1} D_i^{p-1-k}
\Theta(D_i,D_i^{k}(x_i^{p-1}D_j))
=- \Theta(D_i,D_j)=-D_{ij}(x^{\tau}).$$
Therefore, according to Lemma \ref{S-inner}(i), the derivation $H_{\Theta}(D_i)$
is not inner and hence $\Theta\not\in H_{*}^2(S(n),S(n))$ by the Hochschild
exact sequence (\ref{inclusion-sequence}).

On the other hand, we are going to prove that ${\rm Sq}(D_h)\in H_*^2(S(n),S(n))$
(for any $h$) by showing that for any $D_{rs}(x^a)\in S(n)$ the
derivation  $H_{{\rm Sq}(D_h)}(D_{rs}(x^a))$ satisfies the two conditions of
Lemma \ref{S-inner}.

Suppose first, by contradiction, that the first condition of
Lemma \ref{S-inner} is not satisfied for certain indices $i\neq j$, that is
$[H_{{\rm Sq}(D_h)}(D_{rs}(x^a))\cdot x_i^{p-1}D_j]_{D_{ij}(x^{\tau})}\neq 0$.
Then the index $i$ must be equal to $r$ or $s$ and therefore, since we can choose
the index $j\neq i$, we can assume without loss of generality that $(i,j)=(r,s)$.
However, from the definition of the operator $H$, it is straightforward to see that
$$H_{{\rm Sq}(D_h)}(D_{rs}(x^a))\cdot x_r^{p-1}D_s \in
\langle x^{pa-\ep_r-p\ep_s-p\ep_h}D_s, x^{pa-(p+1)\ep_s-p\ep_h}D_r\rangle_F,$$
and this contradicts the hypothesis since the multi-index $pa-\ep_r-p\ep_s-p\ep_h$
cannot be equal to the multi-index $\tau-\ep_r$.

Suppose next, again by contradiction, that the second condition of Lemma \ref{S-inner}
is not satisfied, that is $\sum_k [H_{{\rm Sq}(D_h)}(D_{rs}(x^a))\cdot D_k]_{D_k}\neq 0$.
Then the element $D_{rs}(x^a)$ must have degree $1$ and the formula (\ref{H-morphism})
for $H$ simplifies as
$$H_{{\rm Sq}(D_h)}(D_{rs}(x^a))\cdot D_k ={\rm Sq}(D_h)(D_{rs}(x^a),
{\rm ad} \, D_{rs}(x^a)^{p-1} \cdot D_k).
$$
From this formula it is straightforward to see that
$$H_{{\rm Sq}(D_h)}(D_{rs}(x^a))\cdot D_k \in \langle
D_{rs}(x^{pa-p\ep_h-(p-1)(\ep_r+\ep_s)-\ep_k})\rangle_F.$$
Therefore if $[H_{{\rm Sq}(D_h)}(D_{rs}(x^a))\cdot D_k]_{D_k}\neq 0$, we must have
that $x^a=x_rx_sx_h$ and that $k=r$ or $s$. Now we distinguish two cases, according
to whether $h$ is equal to one of the two indices $r$ and $s$, or not.
If $h\neq r,s$, using the formulas
$$\begin{sis}
& {\rm ad }\, D_{rs}(x_rx_sx_h)^{p-1} \cdot D_r =D_{rs}(x_sx_h^{p-1}), \\
& {\rm ad }\, D_{rs}(x_rx_sx_h)^{p-1} \cdot D_s =- D_{rs}(x_rx_h^{p-1}),
\end{sis}$$
we get a contradiction with the non-vanishing hypothesis because of the following
$$\begin{sis}
& H_{{\rm Sq}(D_h)}(D_{rs}(x_rx_sx_h))\cdot D_r ={\rm Sq}(D_h)(D_{rs}(x_rx_sx_h),
D_{rs}(x_s x_h^{p-1}))=\\
& \hspace{4,1cm} = [D_{rs}(x_rx_s),D_{rs}(x_s)]=-D_r, \\
& H_{{\rm Sq}(D_h)}(D_{rs}(x_rx_sx_h))\cdot D_s ={\rm Sq}(D_h)(D_{rs}(x_rx_sx_h),
-D_{rs}(x_r x_h^{p-1}))=\\
& \hspace{4,1cm} =- [D_{rs}(x_rx_s),D_{rs}(x_r)]=D_s. \\
\end{sis}$$
On the other hand, if $h=r\neq s$, one can prove by induction on
$1\leq t \leq p-1$ that
$$\begin{sis}
& {\rm ad} \, D_{hs}(x_h^2x_s)^t\cdot D_h=\prod_{u=1}^t (u-3) \cdot D_{hs}(x_h^t x_s),\\
& {\rm ad} \, D_{hs}(x_h^2x_s)^{t}\cdot D_s=-t! \, D_{hs}(x_h^{t+1}).\\
\end{sis}$$
Therefore both of the above expressions vanish for $t=p-1$, which contradicts
with our assumption. \qed

\section{The contact algebra}

Let us recall the definition of the restricted contact algebra.
Fix an odd integer $n=2m+1\geq 3$ and a field $F$ of characteristic $p\geq 5$.
For any $j\in\{1,\cdots,2m\}$, we define the sign
$\sigma(j)$ and the conjugate $j'$ of $j$ as follows:
$$\sigma(j)=\begin{sis}
1& \quad \text{ if }1\leq j\leq m,\\
-1& \quad \text{ if } m< j\leq 2m,\\
\end{sis}
\hspace{0,5cm} \text{ and } \hspace{0,5cm}
j'=\begin{sis}
j+m& \quad \text{ if } 1\leq j\leq m,\\
j-m& \quad \text{ if } m<j\leq 2m.
\end{sis}$$
Consider the operator $D_H:A(n)\to W(n)$ defined as
$$D_H(f)=\sum_{j=1}^{2m} \sigma(j)D_j(f)D_{j'}=\sum_{i=1}^m\left[D_i(f)D_{i+m}-D_{i+m}(f)
D_i\right],$$
where, as usual, $D_i:=\frac{\partial}{\partial x_i}\in W(n)$.
We denote by $K'(n)=K'(n,\underline{1})$ the graded Lie algebra over $F$ whose
underlying $F$-vector space is
$A(n)$,
endowed with the grading defined
by ${\rm deg}(x^a)=|a|+a_n-2=\sum_{i=1}^{2m}a_i +2a_n-2$ and with the Lie bracket defined by
\begin{equation*}\label{K-brac}
[x^a,x^b]=D_H(x^a)(x^b)+\left[a_n {\rm deg}(x^b)-b_n {\rm deg}(x^a)\right]
x^{a+b-\epsilon_n}.
\end{equation*}

\begin{defi}
The contact algebra is the derived subalgebra of $K'(n)$:
$$K(n)=K(n;\underline{1})=K'(n)^{(1)}=[K'(n),K'(n)].$$
\end{defi}
Indeed it turns out that
 $$K(n)=
\begin{sis}
&K'(n)& \quad \text{ if } p\not | \: (m+2),\\
&K'(n)_{\neq \tau}& \quad \text{ if } p\: | \:(m+2),\\
\end{sis}$$
where $K'(n)_{\neq \tau}$ is the sub-vector space of $K'(n)$ generated over $F$
by the monomials $x^a$ such that $a\neq \tau:=(p-1,\cdots,p-1)$. Note that
the elements of negative degree are $K(n)_{-2}=\langle 1\rangle_F$ and
$K(n)_{-1}=\langle x_1,\cdots, x_n\rangle_F$.
The $[p]$-map on the base $\{x^a\}$ is given by
$$(x^a)^{[p]}=\begin{sis}
& x^a & \text{ if } x^a=x_ix_{i'} \:\text{ or } \: x^a=x_n, \\
& 0 & \text{ otherwise. }
\end{sis}$$
Observe that our notations for the contact algebra differ slightly from the ones of
\cite[Section 4.5]{FS}, since we have dropped the operator $D_K$ used there.
Therefore, our elements $x^a$ correspond to their elements $D_K(x^a)$.

\hspace{0,4cm}
The first cohomology group of the adjoint representation is equal to (see \cite{CEL}
or \cite[Section 4.8]{FS}):
\begin{equation}
H_{*}^1(K(n),K(n))=H^1(K(n),K(n))=\begin{cases}
0 & \: \text{ if } p\not | \: (m+2),\\
\langle{\rm ad}\:x^{\tau}\rangle_F & \: \text{ if } p \: | \: (m+2).\\
\end{cases}\end{equation}
From this result, we can deduce a criterion saying when a derivation of $K(n)$ is inner in the
case $p | (m+2)$ (in the other case, every derivation is inner).
As usual, for any two elements $E, x^a \in K(n)$, we indicate with $E_{x^a}$
the coefficient of $E$ with respect to the base element $x^a$.

\begin{lem}\label{K-inner}
Suppose that $p$ divides $m+2$. Then a derivation $\gamma:K(n)\to K(n)$ is inner if
and only if $\gamma(1)_{x^{\tau-\ep_n}}=0$.
\end{lem}
\dem
We first prove that the condition is necessary.
Consider an inner derivation ${\rm ad}(x^a)$, with $x^a\in K(n)$.
Then, from the computation
${\rm ad}(x^a)(1)=[x^a,1]=$ $=-2a_n x^{a-\ep_n},$
we deduce that ${\rm ad}(x^a)(1)_{x^{\tau-\ep_n}}= 0$
since $x^{\tau}\not\in K(n)$ by the hypothesis $p | (m+2)$.
However, the condition is also sufficient since ${\rm ad}(x^{\tau})(1)=2x^{\tau-\ep_n}$.
\qed

In \cite[Theorem 1.1]{VIV2}, we prove that the second cohomology group of the
adjoint representation is
\begin{equation}\label{K-ord-def}
H^2(K(n),K(n))=\bigoplus_{i=1}^{2m} \langle {\rm Sq}(x_i) \rangle_F \oplus
\langle {\rm Sq}(1)\rangle_F.
\end{equation}
Using this result and the exact sequence (\ref{inclusion-sequence}), we compute the
second restricted cohomology group.

\dem[of Theorem \ref{K-finaltheorem}]
The theorem is clearly true in the case when $p$ does not divide $m+2$, because
in this case the first cohomology group vanishes. In the case where $p \, | \, (m+2)$,
we are going to show that for any $x^a\in K(n)$ the derivations $H_{{\rm Sq}(1)}(x^a)$
and $H_{{\rm Sq}(x_i)}(x^a)$ satisfy the condition of Lemma \ref{K-inner}.

Consider first the cocycle ${\rm Sq}(1)=2{\rm Sq}(D_n)$. From the definition
(\ref{H-morphism}) of the operator $H$, it is straightforward to check that
$$H_{{\rm Sq}(1)}(x^a)\cdot 1\in \langle x^{pa-2p\ep_n}\rangle_F.$$
Therefore the condition of Lemma \ref{K-inner} is satisfied since the multi-index
$pa-2p\ep_n$ cannot be equal to $\tau-\ep_n$.

Consider next the cocycle ${\rm Sq}(x_i)$.  From the definition
(\ref{H-morphism}) of $H$ together with the fact that ${\rm ad}(x_i)=
\sigma(i)D_{i'}+x_iD_n$, it is straightforward to check that
$$H_{{\rm Sq}(x_i)}(x^a)\cdot 1 \in \langle \sum_{r, s\in \Z}
x^{pa+r\ep_i-(p-r)\ep_{i'}-s\ep_n}\rangle_F.$$
Since $n=2m+1\neq 3$ by the hypothesis $p\, |  \, (m+2)$ and $p\geq 5$,
the multi-index $pa+r\ep_i-(p-r)\ep_{i'}-s\ep_n$ cannot be equal to $\tau-\ep_n$
and therefore the condition of Lemma \ref{K-inner} is satisfied.
\qed

\section{The Hamiltonian algebra}

Let us recall the definition of the restricted Hamiltonian algebra.
Fix an even integer $n=2m\geq 2$ and a field $F$ of
characteristic $p\geq 5$.

\hspace{0,4cm}
We introduce some notations that will be used in this section.
As before, for any $j\in\{1,\cdots,2m\}$, we define the sign
$\sigma(j)$ and the conjugate $j'$ of $j$ as follows:
$$\sigma(j)=\begin{sis}
1& \quad \text{ if }1\leq j\leq m,\\
-1& \quad \text{ if } m< j\leq 2m,\\
\end{sis}
\hspace{0,5cm} \text{ and } \hspace{0,5cm}
j'=\begin{sis}
j+m& \quad \text{ if } 1\leq j\leq m,\\
j-m& \quad \text{ if } m<j\leq 2m.
\end{sis}$$
Given two $n$-tuples of natural numbers $a=(a_1,\cdots, a_n)$ and  $b=(b_1,\cdots,b_n)$,
we say that $a\leq b$ if $a_i\leq b_i$ for every $i$.
We define the degree of $a\in \N^n$ as $|a|=\sum_{i=1}^n a_i$ and the
factorial as $a!=\prod_{i=1}^n a_i!$. For two multi-indices $a, b\in \N^n$ such
that $b\leq a$, we set
$\binom{a}{b}:=\prod_{i=1}^n \binom{a_i}{b_i}=\frac{a!}{b!(a-b)!}.$
Moreover, we define the sign of $a\in \N^{2m}$ as
$\sigma(a)=\prod\sigma(i)^{a_i}$ and the conjugate of $a$ as the multi-index
$\hat{a}$ such that $\hat{a}_i=a_{i'}$ for every $1\leq i\leq 2m$.
We set $\tau:=(p-1,\cdots,p-1)$ (as usual) and
$\0:=(0,\cdots, 0)$.

\hspace{0,4cm}
We denote by $\ti{H(n)}=\ti{H(n;\underline{1})}$ the graded $F$-Lie algebra whose
underlying vector space is $A(n)$, endowed with the grading
defined by ${\rm deg}(x^a)=|a|-2$ and with the Lie bracket defined by
\begin{equation*}
[x^a,x^b]=D_H(x^a)(x^b),
\end{equation*}
where  $D_H: A(n)\to W(n)$ is defined (as before) by
\begin{equation*}D_H(f)=\sum_{j=1}^{2m} \sigma(j)D_j(f)D_{j'}=\sum_{i=1}^m\left[D_i(f)D_{i+m}-
D_{i+m}(f)D_i\right].
\end{equation*}
We denote by $H'(n)=H'(n,\underline{1})$ the quotient of $\ti{H(n)}$ by the central element
$1$.

\begin{defi}
The restricted Hamiltonian algebra is the derived subalgebra of $H'(n)$:
$$H(n)=H(n;\underline{1})=H'(n)^{(1)}=[H'(n),H'(n)].$$
\end{defi}
Observe that $H(n)$ has $F$-dimension $p^n-2$, with a base given by the elements $\{x^a\}$
such that $x^a\neq 1$ and $x^a\neq x^{\tau}$. The elements of negative degree
are $H(n)_{-1}=\langle x_1,\cdots,x_n \rangle_F$. On the elements of the base,
the $[p]$-map is given by
\begin{equation*}
(x^a)^{[p]}=\begin{sis}
& x^a & \text{ if } x^a=x_ix_{i'}, \\
& 0 &\text {otherwise. }
\end{sis}
\end{equation*}
Note that our notations for the Hamiltonian algebra differ slightly from the ones of
\cite[Section 4.4]{FS}, since for simplicity we have dropped the operator $D_H$
used there. Therefore, our elements $x^a$ correspond to their elements $D_H(x^a)$.

\hspace{0,4cm}
The first cohomology group of the adjoint representation is given by (see \cite{CEL}
or \cite[Section 4.8]{FS}):
\begin{equation}
H_{*}^1(H(n),H(n))=H^1(H(n),H(n))=\langle {\rm ad} \, x^{\tau}\rangle_F
\bigoplus_{i=1}^n \langle x_i^{p-1}D_{i'}\rangle_F
\oplus \langle {\rm deg}\rangle_F,
\end{equation}
where $x_i^{p-1}D_{i'}$ is the derivation which sends $x^a\in H(n)$ into $x_i^{p-1}D_{i'}(x^a)
\in H(n)$ and ${\rm deg}$ is the operator degree defined by ${\rm deg}(x^a)={\rm deg}(x^a)\:
x^a$.

\hspace{0,4cm}
From this result, we can deduce a criterion saying when a derivation of $H(n)$ is inner.
As usual, given two elements $E, x^{a} \in H(n)$, we denote by $E_{x^a}$ the coefficient
of $E$ with respect to the base element $x^a$.

\begin{lem}\label{H-inner}
A derivation $\gamma:H(n)\to H(n)$ is inner if and only if it satisfies the three following
conditions:
\begin{itemize}
\item[(i)] There exists an index $i$ such that $\gamma(x_i)_{x^{\tau-\ep_{i'}}}=0$.
\item[(ii)] For any index $i$, it holds that
$\gamma(x^{\tau-(p-1)\ep_i})_{x^{\tau-\ep_{i'}}}=0$.
\item[(iii)] There exists an index $i$ such that $\gamma(x_i)_{x_i}+\gamma(x_{i'})_{x_{i'}}=0$.
\end{itemize}
\end{lem}
\dem
We first prove that the three conditions are necessary. Consider an inner derivation
${\rm ad}(D)$, with $D\in H(n)$. Write $D$ as linear combination of the base elements $x^a$.
The element ${\rm ad}(x^a)(x_i)=[x^a,x_i]=-\sigma(i)a_{i'}x^{a-\ep_{i'}}$ cannot belong
to $\langle x^{\tau-\ep_{i'}}\rangle_F$ since $x^{\tau}\not\in H(n)$. Therefore the condition
(i) is verified. Consider now the element
$$[x^a,x^{\tau-(p-1)\ep_i}]=-\sigma(i)a_i x^{a+\tau-p\ep_i-\ep_{i'}}+
\sum_{j\neq i, i'} \sigma(j)[a_{j'}-a_j]x^{a+\tau-(p-1)\ep_i-\ep_j-\ep_{j'}}.
$$
We have that $x^{a+\tau-p\ep_i-\ep_{i'}}\not\in \langle x^{\tau-\ep_{i'}}\rangle_F$
since $a_i\leq p-1$ and  $x^{a+\tau-(p-1)\ep_i-\ep_j-\ep_{j'}}\not \in
\langle x^{\tau-\ep_{i'}}\rangle_F$  since $a_{i'}\geq 0$. Therefore condition (ii)
is verified.
Finally, to prove condition (iii), we write $D=\sum_{i=1}^m a_i x_ix_{i'} +E$ with
$E_{x_ix_{i'}}=0$ for every $i$. For any index $i$, we have that
$$D(x_i)_{x_i}+D(x_{i'})_{x_{i'}}=-\sigma(i)a_i+\sigma(i)a_i=0,$$
and therefore also condition (iii) is verified.
However, the three conditions are also sufficient since
$$\begin{sis}
&{\rm ad}(x^{\tau})(x_i)=[x^{\tau},x_i]=\sigma(i)x^{\tau-\ep_{i'}},\\
&(x_i^{p-1}D_{i'})(x^{\tau-(p-1)\ep_i})=-x^{\tau-\ep_{i'}},\\
& {\rm deg}(x_i)_{x_i}+{\rm deg}(x_{i'})_{x_{i'}}=-1-1=-2.\qed
\end{sis}$$

In \cite[Theorem 1.2]{VIV1}, we prove that the second cohomology group of the
adjoint representation is
\begin{equation}\label{H-ord-def}
H^2(H(n),H(n))=\begin{sis}
&\bigoplus_{i=1}^n \langle {\rm Sq}(x_i)\rangle_F
\bigoplus_{\stackrel{i<j}{j\neq i'}} \langle \Pi_{ij} \rangle_F \bigoplus_{i=1}^m
\langle \Pi_i \rangle_F \oplus \langle \Phi\rangle_F & \text{if }n\geq 4,\\
&\bigoplus_{i=1}^2 \langle {\rm Sq}(x_i)\rangle_F
\oplus \langle \Phi\rangle_F &\text{ if }n=2,\\
\end{sis}
\end{equation}
where the above cocycles are defined (and vanish elsewhere) by
$$\begin{sis}
& \Pi_{ij}(x^a, x^b)=x_{i'}^{p-1}x_{j'}^{p-1}[D_i(x^a)D_j(x^b)-D_i(x^b)D_j(x^a)]
\hspace{0,3cm} \text{ for } j\neq i,i',  \\
& \Pi_i(x_ix^a,x_{i'}x^b)=x^{a+b+(p-1)(\ep_i+\ep_{i'})} \hspace{0,5cm}
\text{ if } a+b\neq \tau-(p-1)(\ep_i+\ep_{i'}),\\
& \Pi_i(x_k,x^{\tau-(p-1)(\ep_i+\ep_{i'})})=-\sigma(k)x^{\tau-\ep_{k'}} \hspace{1,5cm}
\text{ for any }\hspace{0,3cm} 1\leq k\leq n,\\
& \Phi(x^a,x^b)=\sum_{\stackrel{\0\leq \delta\leq a, \widehat{b}}{|\delta|=3}}
\binom{a}{\delta}\binom{b}{\widehat{\delta}}\sigma(\delta)\:\delta!\:
x^{a+b-\delta -\widehat{\delta}}.\\
\end{sis}$$

Using the above result (\ref{H-ord-def}) and the Hochschild exact sequence
(\ref{inclusion-sequence}), we can compute the second restricted cohomology group.

\dem[of the Theorem \ref{H-finaltheorem}]
For $n\geq 4$, the cocycles $\Pi_{ij}$ and $\Pi_i$ do not belong to
$H_{*}^2(H(n), $ $H(n))$. To unify the notation we set (only during this proof)
$\Pi_{ii'}=\Pi_i$ for $1\leq i\leq m$. Then, using that $x_r^{[p]}=0$,
we compute (for $i<j$):
$$\begin{sis}
& H_{\Pi_{ij}}(x_i)\cdot x^{\tau-(p-1)\ep_{j'}}=\Pi_{ij}(x_i,
({\rm ad} x_i)^{p-1}(x^{\tau-(p-1)\ep_{j'}}))=\sigma(i)x^{\tau-\ep_j},\\
& H_{\Pi_{ij}}(x_j)\cdot x^{\tau-(p-1)\ep_{i'}}=\Pi_{ij}(x_j,
({\rm ad} x_j)^{p-1}(x^{\tau-(p-1)\ep_{i'}}))=-\sigma(j)x^{\tau-\ep_i},\\
& H_{\Pi_{ij}}(x_k)\cdot x^{\tau-(p-1)\ep_h}=0 \hspace{3cm} \text{ if } k\neq i, j
\text{ or } h\neq i', j'.
\end{sis}$$
Therefore according to Lemma \ref{H-inner}(ii), no linear combination of $\Pi_{ij}$
can be in the kernel of the map $H$ and therefore in $H_{*}^2(H(n),H(n))$.

\hspace{0,4cm}
Next we prove that ${\rm Sq}(x_h)=\sigma(h){\rm Sq}(D_{h'})
\in H_*^2(H(n),H(n))$ by showing
that for any $x^a\in H(n)$ the derivation $H_{{\rm Sq}(x_h)}(x^a)$ satisfies the
conditions of Lemma \ref{H-inner}. Suppose,
by contradiction, that the condition (ii) of Lemma \ref{H-inner}
is not satisfied for some index $i$. Then for degree reasons we must have
$$p-2={\rm deg}(x_i^{p-1}D_{i'})={\rm deg}({\rm Sq}(x_h))+p \,{\rm deg}(x^a)=-p+
p\, {\rm deg}(x^a),$$
which is impossible since $p\neq 2$. Suppose now, by contradiction, that the
condition (iii) of Lemma \ref{H-inner} is not satisfied and in particular that
\begin{equation*}
[H_{{\rm Sq}(x_h)}(x^a)\cdot x_h]_{x_h}+[H_{{\rm Sq}(x_h)}(x^a)\cdot x_{h'}]_{x_{h'}}
\neq 0.\tag{*}
\end{equation*}
From the definition of $H$, it is straightforward to see that
\begin{equation*}\tag{**}\begin{sis}
&  H_{{\rm Sq}(x_h)}(x^a)\cdot x_h\in \langle x^{pa-2p\ep_{h'}-(p-1)\ep_h}\rangle_F,\\
&  H_{{\rm Sq}(x_h)}(x^a)\cdot x_{h'}\in
\langle x^{pa-(2p-1)\ep_{h'}-p \ep_h}\rangle_F.\\
\end{sis}\end{equation*}
From the hypothesis (*), it follows that $x^a=x_{h'}^2x_h$ and therefore the formula
(\ref{H-morphism}) simplifies (for any index $k$) as follows
$$H_{{\rm Sq}(x_h)}(x_{h'}^2 x_h)\cdot x_k= {\rm Sq}(x_h)(x_{h'}^2x_h,({\rm ad}\,
x_{h'}^2x_h)^{p-1} x_k).$$
By induction on $1\leq r\leq p-1$, one can verify that
$$\begin{sis}
& ({\rm ad}\, x_{h'}^2x_h)^r x_{h'}=\sigma(h')^r(-1)(-2)\cdots (-r) \, x_{h'}^{r+1},\\
& ({\rm ad}\, x_{h'}^2x_h)^r x_h=\sigma(h')^r 2 (2-1)\cdots (2-(r-1))\, x_{h'}^r x_h.
\end{sis}$$
In particular we have that $({\rm ad}\, x_{h'}^2x_h)^{p-1} x_{h'}=
({\rm ad}\, x_{h'}^2x_h)^{p-1} x_h=0$ and, substituting in the above expression, we get
$H_{{\rm Sq}(x_h)}(x_{h'}^2 x_h)\cdot x_{h'}=
H_{{\rm Sq}(x_h)}(x_{h'}^2 x_h)\cdot x_h$ contradicting the hypothesis (*).
Finally, by using the first equation of (**) and the fact that
the multi-index $pa-2p\ep_{h'}-(p-1)\ep_h$ cannot be equal to $\tau-\ep_{h'}$,
the condition (i) of Lemma \ref{H-inner} is satisfied for $i=h$.

\hspace{0,4cm}
Finally, we prove that $\Phi\in H_*^2(H(n),H(n))$ by showing that for any
$x^a\in H(n)$ the derivation $H_{\Phi}(x^a)$
satisfies the conditions of Lemma \ref{H-inner}. Suppose, by contradiction, that the
condition (ii) of Lemma \ref{H-inner} is not satisfied for some index $i$.
Then for degree reasons, we must have
$$p-2={\rm deg}(x_i^{p-1}D_{i'})={\rm deg}(\Phi)+p\, {\rm deg}(x^a)=-4+p\, {\rm deg}(x^a),$$
a contradiction. Analogously, if the condition (iii) of the Lemma is
not satisfied by $H_{\Phi}(x^a)$ for some index $i$, then we get a contradiction
by looking at the degree
$$ 0={\rm deg}(\Phi)+p\, {\rm deg}(x^a)=-4 +p\, {\rm deg}(x^a).$$
Finally, suppose that the condition (i) of Lemma \ref{H-inner} is not satisfied for
some index $i$, that is $[H_{\Phi}(x^a)\cdot x_i]_{x^{\tau-\ep_{i'}}}\neq 0$.
Then, by looking at the degree, we get that
$${\rm deg}(x^a)=2m-2\frac{m-1}{p}>0.$$
In particular we have that $p \, |\, (m-1)$, from which we deduce that either $m=1$
or $m\geq p+1\geq 6$. Suppose first that $m\neq 1$. Then, from the formula
(\ref{H-morphism}) and  using that $(x^a)^{[p]}=0$, we deduce that
$$H_{\Phi}(x^a)\cdot x_i \in \sum_{|\delta|=3} \sum_{k=0}^{p-1}
\langle ({\rm ad } x^a)^{p-1-k}(x^{(k+1)a-(k-1)\ep_i-k\ep_{i'}-\delta-\widehat{\delta}})
\rangle_F.$$
Fix a multi-index $\delta$ appearing in the above summation and choose an index $j\neq i, i'$
such $\delta_j=\delta_{j'}=0$ (this is possible since $|\delta|=3$ and $n=2m\geq 12$).
Then the $j$-th coefficient of every monomial appearing in the expression
$$\sum_{k=0}^{p-1}({\rm ad } x^a)^{p-1-k}(x^{(k+1)a-(k-1)\ep_i-k\ep_{i'}-
\delta-\widehat{\delta}})$$
is $p\, a_j\neq p-1$. Therefore the monomial $x^{\tau-\ep_{i'}}$ cannot appear in the above
expression and, repeating the same argument for every multi-index $\delta$ as before, we get
that $[H_{\Phi}(x^a)\cdot x_i]_{x^{\tau-\ep_{i'}}}=0$, a contradiction. In the remaining case
$m=1$, we have that
$$H_{\Phi}(x^a)\cdot x_i \in \sum_{|\delta|=3}
\langle x^{pa -(p-2)\ep_i-(p-1)\ep_{i'}-\delta-\hat{\delta}}  \rangle_F.$$
From this and the hypothesis $[H_{\Phi}(x^a)\cdot x_i]_{x^{\tau-\ep_{i'}}}\neq 0$,
we deduce that $x^a=x_i^2x_{i'}^2$.
Using the straightforward formulas
$$\begin{sis}
& ({\rm ad}\, x_i^2x_{i'}^2)^k(x_i) = [-2\sigma(i)]^k x_i^{k+1}x_{i'}^k,\\
& \Phi(x_i^2x_{i'}^2,x_i^{k+1}x_{i'}^k)=2\sigma(i)(k+1)k \, x_i^k x_{i'}^{k-1},\\
& ({\rm ad}\, x_i^2x_{i'}^2)^{p-1-k}(x_i^kx_{i'}^{k-1}) = [-2\sigma(i)]^{p-1-k}
 x_i^{p-1}x_{i'}^{p-2},
\end{sis}$$
we get that
$H_{\Phi}(x^a)\cdot x_i=2\sigma(i) \sum_{k=0}^{p-1} [(k+1)k] x^{\tau-\ep_{i'}}=0,$
since $\sum_{k=0}^{p-1} k \equiv 0 \mod p$ for $p\geq 3$ and
$\sum_{k=0}^{p-1} k^2\equiv 0 \mod p$ for $p\geq 5$.
This contradiction finishes the proof.
\qed

\section{The Melikian algebra}

Let us recall the definition of the restricted Melikian algebra, following
\cite[Section 4.3]{STR}.

\hspace{0,4cm}
Let $F$ be a field of characteristic $p=5$.
Consider $W(2)={\rm Der}_F A(2)={\rm Der}_F F[x_1,x_2]/(x_1^p,x_2^p)$,
the restricted Witt-Jacobson Lie algebra of rank $2$.
Let $\ti{W(2)}$ be a copy of $W(2)$
and for an element $D\in W(2)$ we indicate with $\widetilde{D}$ the corresponding element
inside $\widetilde{W(2)}$. The restricted Melikian algebra $M:=M(1,1)$ is defined as
$$M=A(2)\oplus W(2)\oplus \widetilde{W(2)},$$
with Lie bracket defined by the following rules (for all $D, E\in W(2)$ and
$f, g\in A(2)$):
$$\begin{sis}
&[D,\widetilde{E}]:=\widetilde{[D,E]}+2\,{\rm div}(D)\widetilde{E},\\
&[D,f]:=D(f)-2\, {\rm div}(D)f,\\
&[f_1\widetilde{D_1}+f_2\widetilde{D_2},g_1\widetilde{D_1}+g_2\widetilde{D_2}]:=
f_1g_2-f_2g_1,\\
&[f,\widetilde{E}]:=f E,\\
&[f,g]:=2\,(gD_2(f)-fD_2(g))\widetilde{D_1}+2\,(fD_1(g)-gD_1(f))\widetilde{D_2},\\
\end{sis}$$
where ${\rm div}(f_1D_1+f_2D_2):=D_1(f_1)+D_2(f_2)\in A(2)$.
The Melikian algebra $M$ has
a $\Z$-grading given by (for all $D, E\in W(2)$ and
$f\in A(2)$):
$$\begin{sis}
&{\rm deg}_M(D):=3\,{\rm deg}(D),\\
&{\rm deg}_M(\widetilde{E}):=3\,{\rm deg}(E)+2,\\
&{\rm deg}_M(f):=3\,{\rm deg}(f)-2.\\
\end{sis}$$
In particular the elements of negative degree are
$$M_{-3}=\langle D_1, D_2 \rangle_F, \hspace{1cm} M_{-2}=\langle 1\rangle_F ,
\hspace{1cm} M_{-1}=\langle \widetilde{D_1}, \widetilde{D_2} \rangle_F.$$
The $[p]$-map is defined on an element $X\in M$ by
$$X^{[p]}=\begin{sis}
& X & \text{ if } X=x_1D_1 \text{ or } x_2D_2 , \\
& 0 & \text{ otherwise.}
\end{sis}$$
It is known that every derivation of $M$ is inner (see \cite[Section 4.3]{STR}), that is
\begin{equation}
H_{*}^1(M,M)=H^1(M,M)=0.
\end{equation}
Therefore, from the Hochschild exact sequence (\ref{inclusion-sequence}), we deduce that
$H_*^2(M,$ $M)$ $=H^2(M,M)$. The Theorem \ref{M-finaltheorem} follows from \cite[Theorem 1.1]{VIV3}:
$$H^2(M,M)=\langle {\rm Sq}(1)\rangle_F \bigoplus_{i=1}^2 \langle {\rm Sq}(D_i)
\rangle_F \bigoplus_{i=1}^2 \langle {\rm Sq}(\ti{D_i})\rangle_F.$$

\section*{Acknowledgements}
The author was supported by the grant FCT-Ci\^encia2008 from CMUC (University of Coimbra)
and by the FCT project \textit{Espa\c cos de Moduli em Geometria Alg\'ebrica} (PTDC/MAT/111332/2009).

\end{document}